\documentclass[letterpaper, 10 pt, conference, onecolumn]{ieeeconf}

\IEEEoverridecommandlockouts

\usepackage{cite}
\usepackage{graphics}
\usepackage{epsfig}
\usepackage{ntheorem,amsmath}
\usepackage{amssymb}
\usepackage{mathtools}
\usepackage{commath}
\usepackage{graphicx,bm}
\usepackage{picture}
\usepackage{xcolor}
\usepackage{float}
\usepackage{verbatimbox}
\usepackage{booktabs}
\usepackage{listings}
\usepackage{tcolorbox}
\usepackage{url}
\usepackage{algorithm}
\usepackage{algpseudocode}
\usepackage{textcomp}
\usepackage{grffile}
\usepackage{times}
\usepackage{hyperref}

\makeatletter
\renewtheoremstyle{plain}
{\item{\theorem@headerfont ##1\ ##2\theorem@separator}~}
{\item{\theorem@headerfont ##1\ ##2\ (##3)\theorem@separator}~}
\theoremseparator{.}
\makeatother

\newtheorem{theorem}{Theorem}

\newtheorem{assumption}{Assumption}

\newtheorem{lemma}{Lemma}
\newtheorem{remark}{Remark}
\newtheorem{definition}{Definition}

\DeclareMathOperator{\col}{col}
\DeclareMathOperator{\diag}{diag}
\DeclareMathOperator{\rank}{rank}

\newcommand{\N}{{\mathbb{N}}}
\newcommand{\R}{{\mathbb{R}}}
\newcommand{\C}{{\mathbb{C}}}
\newcommand{\Ts}{T_{\textup{s}}}

\renewcommand{\QED}{\QEDopen}

\title{\LARGE \bf Derivative-Free Data-Driven Control of Continuous-Time\\
Linear Time-Invariant Systems
}

\author{Alessandro Bosso$^{1}$, Marco Borghesi$^{1}$, Andrea Iannelli$^{2}$, Giuseppe Notarstefano$^{1}$, and Andrew R. Teel$^{3}$%
\thanks{$^{1}$A. Bosso, M. Borghesi, and G. Notarstefano are with the Department of Electrical, Electronic, and Information Engineering, University of Bologna, Italy. Email: {\tt\small \{alessandro.bosso, m.borghesi, giuseppe.notarstefano\}@unibo.it}}
\thanks{$^{2}$A. Iannelli is with the Institute for System Theory and Automatic Control, University of Stuttgart, Germany. Email: {\tt\small andrea.iannelli@ist.uni-stuttgart.de}
}
\thanks{$^{3}$A. R. Teel is with the Department of Electrical and Computer Engineering, University of California, Santa Barbara, CA, USA. Email: {\tt\small teel@ece.ucsb.edu}}
\thanks{The research leading to these results has received funding from the European Union's Horizon Europe research and innovation program under the Marie Sk{\l}odowska-Curie Grant Agreement No. 101104404 - \mbox{IMPACT4Mech}.}
}

\begin{document}

    \maketitle
	\thispagestyle{empty}
	\pagestyle{empty}
	
	\begin{abstract}
        This paper develops a data-driven stabilization method for continuous-time linear time-invariant systems with theoretical guarantees and no need for signal derivatives.
        The framework, based on linear matrix inequalities (LMIs), is illustrated in the state-feedback and single-input single-output output-feedback scenarios.
        Similar to discrete-time approaches, we rely solely on input and state/output measurements.
        To avoid differentiation, we employ low-pass filters of the available signals that, rather than approximating the derivatives, reconstruct a non-minimal realization of the plant.
        With access to the filter states and their derivatives, we can solve LMIs derived from sample batches of the available signals to compute a dynamic controller that stabilizes the plant.
        The effectiveness of the framework is showcased through numerical examples.
	\end{abstract}
    
    \section{Introduction}
        Over the past decades, the paradigm of data-driven learning has gained increasing attention in control theory.
        This interest can be traced back to areas such as system identification \cite{ljung1999system} and adaptive control \cite{ioannou2012robust}, while the recent trends mainly focus on reinforcement learning \cite{sutton2018reinforcement}.
        A common theme across these approaches is the shift from relying on precise models to leveraging the information contained in the collected data.
        Recently, the dominant paradigm in data-driven control has become to compute controllers directly from data using linear matrix inequalities (LMIs) or other optimization problems, without even requiring an intermediate identification step \cite{de2019formulas}.
        In this paper, we focus on LMI-based methods for the stabilization of continuous-time linear time-invariant systems.

        Fundamental contributions to data-driven control of discrete-time systems include \cite{de2019formulas} and \cite{van2020noisy}, which introduced two distinct data-based LMI formulations for state-feedback stabilization.
        These methodologies have also been used to address the stabilization of bilinear systems \cite{bisoffi2020data}, linear time-varying systems \cite{nortmann2020data}, and the linear quadratic regulator problem \cite{de2019formulas}, also accounting for the effects of noise \cite{de2021low, dorfler2023certainty}.
        The integration of partial model knowledge into these approaches was explored in \cite{berberich2022combining}.
        In this context, necessary and sufficient conditions for data informativity have been thoroughly investigated \cite{van2020data}. 
        
        In the continuous-time scenario, the discrete-time state-feedback stabilization paradigm can be recovered via suitable sampling techniques.
        However, this comes at the cost of requiring samples of the state derivatives \cite{de2019formulas}, causing robustness issues in the presence of noise.
        In \cite{berberich2021data}, LMIs inspired by \cite{van2020noisy} were derived for the design of a stabilizing gain with non-periodic sampling and noisy state-derivative estimates.
        Similarly, derivative estimates were employed in \cite{miller2022data}, which proposed quadratic matrix inequalities for stabilizing linear parameter-varying systems in both discrete and continuous time.
        Recent contributions in continuous time include the study of how sampling impacts data informativity \cite{eising2024sampling} and the stabilization of switched and constrained systems \cite{bianchi2025data}.

        To avoid differentiation in the state-feedback scenario, \cite{de2023event} proposed to construct datasets based on integrals and temporal differences of the available signals.
        To the best of the authors' knowledge, no other strategy in the continuous-time literature completely removes the need for state derivatives.
        Furthermore, no output-feedback approach has been developed thus far, where the sensitivity to noise is even more pronounced due to the need for multiple differentiations.

        This paper presents a data-driven stabilization framework for continuous-time systems that eliminates the need for signal derivatives.
        The approach is first demonstrated for the state-feedback scenario and then extended with minor modifications to the single-input single-output (SISO) output-feedback case.
        Instead of using derivative approximations or methods based on integrals and temporal differences, we define a non-minimal realization of the plant, inspired by adaptive observer design \cite[Ch. 4]{narendra2012stable} and presented here in a state-space setting.
        Specifically, we process the input and state/output signals with low-pass filters that are shown to converge exponentially to an augmented system representation.
        Since both the state and the derivative of the filters are accessible, we employ LMIs similar to those in \cite{de2019formulas} to compute the gains of a dynamic, filter-based, stabilizing controller.
        Feasibility of the LMIs is ensured under suitable excitation conditions, and closed-loop stability is guaranteed regardless of the initial filter transient.
        Numerical examples validate the effectiveness of the approach.

        The paper is organized as follows.
        In Section \ref{sec:problem}, we state the design problem and introduce LMI-based data-driven control.
        In Sections \ref{sec:state_feedback} and \ref{sec:output_feedback}, we provide the algorithms for the state-feedback and the SISO scenarios.
        In Section \ref{sec:numerical_example}, we show the numerical results.
        Finally, Section \ref{sec:conclusion} concludes the paper.

        \emph{Notation:} We use $\N$, $\R$, and $\C$ to denote the sets of natural, real and complex numbers.
        The identity of dimension $j \in \N$ is denoted with $I_j$.
        Given a symmetric matrix $M = M^\top$, $M \succ 0$ (resp. $M \succeq 0$) denotes that it is positive definite (resp. positive semidefinite).
        Similarly, $\prec 0$ and $\preceq 0$ are used for negative definite and semidefinite matrices.
    
    \section{Problem Statement and Preliminaries}\label{sec:problem}
        Although the next sections will deal with both state and output feedback, for convenience, we illustrate the problem for the state-feedback scenario.
        Consider a linear time-invariant system of the form
        \begin{equation}\label{eq:plant_state}
            \dot{x} = Ax + Bu,
        \end{equation}
        where $x\in \R^n$ is the state, $u \in \R^m$ is the control input, and $A$ and $B$ are unknown matrices of appropriate dimensions.
        For a given initial condition $x(0)$ and some input sequence $u(t)$, suppose that the resulting input-state trajectory of \eqref{eq:plant_state} has been collected over an interval $[0, T]$, with $T > 0$.
        More specifically, suppose that the continuous-time dataset
        \begin{equation}\label{eq:dataset}
            (u(t), x(t)), \qquad \forall t \in [0, T],
        \end{equation}
        is available.
        We are interested in finding an algorithm that uses \eqref{eq:dataset} to compute a stabilizing controller for \eqref{eq:plant_state}, without any prior knowledge of $A$ and $B$.
        
        We present some preliminary notions related to the existing approaches in the literature.
        To recover the results of data-driven stabilization of discrete-time systems, algorithms developed in a continuous-time setting are based on collecting a finite batch of data of $u$, $x$, and $\dot{x}$ with a suitable sampling mechanism \cite{de2019formulas, berberich2021data, miller2022data}.
        Given a fixed sampling time $\Ts \coloneqq T/N$, with $N \in \mathbb{N}$, the following batch is obtained:
        \begin{equation}\label{eq:batches}
            \begin{split}
                U &\coloneqq \begin{bmatrix}
                    u(0) & u(\Ts) & \cdots & u((N - 1)\Ts)
                \end{bmatrix} \in \R^{m \times N}\\
                X &\coloneqq \begin{bmatrix}
                    x(0) & x(\Ts) & \cdots & x((N - 1)\Ts)
                \end{bmatrix} \in \R^{n \times N}\\
                \dot{X} &\coloneqq \begin{bmatrix}
                    \dot{x}(0) & \dot{x}(\Ts) & \cdots & \dot{x}((N - 1)\Ts)
                \end{bmatrix} \in \R^{n \times N}.
            \end{split}
        \end{equation}
        We introduce a key definition used in this paper.
        \begin{definition}\label{def:E}
            A data batch of the form \eqref{eq:batches} is exciting if
            \begin{equation}\label{eq:rank_condition}
            \rank \begin{bmatrix} 
                X \\ U
            \end{bmatrix} = n + m.
            \end{equation}
        \end{definition}
        In the scenario where the data \eqref{eq:batches} are exciting, they can be used to construct a stabilizing feedback law for system \eqref{eq:plant_state} of the form $u = K x$.
        In particular, it is possible to make the closed-loop system matrix $A + BK$ Hurwitz by choosing
        \begin{equation}\label{eq:K}
            K = UQ (XQ)^{-1},
        \end{equation}
        where $Q \in \R^{N \times n}$ is any solution of the following LMI:
        \begin{equation}\label{eq:LMI}
            \begin{cases}
                \dot{X}Q + Q^\top \dot{X}^\top \prec 0\\
                XQ = Q^\top X^\top \succ 0.
            \end{cases}
        \end{equation}
        This result follows mutatis mutandis from the discrete-time case; see \cite[Thms. 2 and 3]{de2019formulas}.
        
        We remark that, in the discrete-time scenario, only the data $X$ and $U$ are needed to compute $K$ \cite{de2019formulas}.
        Instead, the continuous-time framework requires $\dot{X}$, which cannot be reliably inferred from \eqref{eq:dataset} if the data are corrupted by noise.
        Also, even in a noise-free scenario, approximation via finite differences leads to persistent errors in the dataset.
        
        The direct data-driven control framework proposed in this paper addresses the above issue both in the state-feedback scenario of system \eqref{eq:plant_state} and in the case of output feedback for SISO systems, where a continuous-time algorithm corresponding to one provided in \cite{de2019formulas} would involve also higher order derivatives of the input and the measured output. 

    \section{Data-Driven Control from Input-State Data}\label{sec:state_feedback}        
        In this section, we are interested in designing a stabilizing controller for system \eqref{eq:plant_state} under the following assumption.
        \begin{assumption}\label{hyp:ctrb}
            The pair $(A, B)$ is controllable.
        \end{assumption}
        To avoid the challenge of having to measure $\dot{x}$, we propose a strategy that involves designing a filter of $x$ and $u$.
        This filter, whose state and state derivative are accessible, is not used to approximate $\dot{x}$ but to reconstruct a non-minimal realization of the plant \eqref{eq:plant_state}.
        Thus, our approach avoids the robustness issues originating from the computation of derivatives from noisy data.
        We begin by presenting the non-minimal realization of system \eqref{eq:plant_state}.
            
        \subsection{Non-Minimal System Realization}\label{sec:filter_state}
            Consider the following dynamical system, having input $u$ and output $\xi \in \R^n$:
            \begin{equation}\label{eq:filter_state_compact}
                \begin{split}
                    \dot{\zeta} &= \begin{bmatrix}
                        A & B\\
                        0 & -\lambda I_m
                    \end{bmatrix} \zeta + \begin{bmatrix}
                        0 \\ \gamma I_m
                    \end{bmatrix} u
                    \\
                    \xi &= \frac{1}{\gamma}\begin{bmatrix}
                        A + \lambda I_n & B
                    \end{bmatrix} \zeta,
                \end{split}
            \end{equation}
            where $\zeta \in \R^{n+m}$ is the state and $\lambda$ and $\gamma$, with $\gamma \neq 0$, are constant scalar tuning gains.
            System \eqref{eq:filter_state_compact} is obtained by compactly rewriting the following system:
            \begin{equation}\label{eq:filter_state}
                \begin{split}
                    \dot{\zeta} &= -\lambda \zeta + \gamma \begin{bmatrix}
                        \xi \\ u
                    \end{bmatrix}\\
                    \xi &= \frac{1}{\gamma}\begin{bmatrix}
                        A + \lambda I_n & B
                    \end{bmatrix} \zeta,
                \end{split}
            \end{equation}
            which, for $\lambda > 0$, acts as a low-pass filter of $\xi$ and $u$.
            The relationship between systems \eqref{eq:plant_state} and \eqref{eq:filter_state_compact} is provided in the next lemma.
            \begin{lemma}\label{lemma:i_o_equivalence_state}  
                Under Assumption \ref{hyp:ctrb}, for all $\lambda$ and all $\gamma \neq 0$, the observable and controllable subsystem of \eqref{eq:filter_state_compact} obeys the same dynamics of \eqref{eq:plant_state}, with state $\xi \in \R^n$.
            \end{lemma}
            In other words, Lemma \ref{lemma:i_o_equivalence_state} states that all input-state trajectories of \eqref{eq:plant_state} are input-output trajectories of \eqref{eq:filter_state_compact}, and vice-versa.
            
            It is also useful to recognize the structural property of system~\eqref{eq:filter_state_compact} given in Lemma \ref{lemma:ctrb_state}.
            \begin{lemma}\label{lemma:ctrb_state}
                Under Assumption \ref{hyp:ctrb}, for all $\lambda$ and all $\gamma \neq 0$, the pair
                \begin{equation}
                    \left(\begin{bmatrix}
                        A & B\\
                        0 &  -\lambda I_m
                    \end{bmatrix}, 
                    \begin{bmatrix}
                        0\\
                        \gamma I_m
                    \end{bmatrix}\right)
                \end{equation}
                is controllable.
            \end{lemma}
            
        \subsection{Controller Design}
            \begin{algorithm}[t!]
        		\caption{Controller Design from Input-State Data}\label{alg:input_state_data}
        		\begin{algorithmic}
        			\State \hspace{-0.47cm} \textbf{Initialization} 
        			\State \emph{Dataset:}
                    \begin{equation}\label{eq:dataset_state}
                        (u(t), x(t)), \qquad \forall t \in [0, T].
                    \end{equation}
                    \State \emph{Tuning:} $\lambda > 0$, $\gamma \neq 0$, $\Ts > 0$.
        			\State \hspace{-0.47cm} \textbf{Data Batches Construction} 
                    \State \emph{Filter of the data:}
                    \begin{equation}\label{eq:filter_state_obs}
                        \dot{\hat{\zeta}}(t) = -\lambda \hat{\zeta}(t) + \gamma \begin{bmatrix}
                            x(t) \\ u(t)
                        \end{bmatrix},\quad \forall t \in [0, T].
                    \end{equation}
                    Initialization: $\hat{\zeta}(0) = 0$.
                    \State \emph{Sampled data batches:}
        			\begin{equation}\label{eq:batches_state}
                        \begin{split}
                            U\! &\coloneqq\! \begin{bmatrix}
                            u(0) & \qquad u(\Ts) & \cdots & \quad u((N - 1)\Ts)
                            \end{bmatrix} \in \R^{m \times N}\\
                            Z\! &\coloneqq\! \begin{bmatrix}
                            \hat{\zeta}(0) & \qquad \hat{\zeta}(\Ts) & \cdots & \quad \hat{\zeta}((N - 1)\Ts)
                            \end{bmatrix} \in \R^{(n + m)\times N}\\
                            \dot{Z}\! &\coloneqq\! \begin{bmatrix}
                            \dot{\hat{\zeta}}(0) & \qquad \dot{\hat{\zeta}}(\Ts) & \cdots & \quad \dot{\hat{\zeta}}((N - 1)\Ts)
                            \end{bmatrix} \in \R^{(n+m) \times N}\\
                            E\! & \coloneqq\! \begin{bmatrix}
                                x(0) & e^{-\lambda\Ts}x(0) & \cdots & e^{-\lambda (N - 1) \Ts}x(0)
                            \end{bmatrix} \in \R^{n \times N}.
                        \end{split}
                    \end{equation}
                    \State \hspace{-0.47cm} \textbf{Stabilizing Gain Computation}
                    \State \emph{LMI:} find $Q \in \R^{N \times (n + m)}$ such that:
                    \begin{equation}\label{eq:LMI_state}
                        \begin{cases}
                            \left(\!\dot{Z} - \begin{bsmallmatrix}
                                \gamma I_n \\ 0
                            \end{bsmallmatrix}E\right)\!Q + Q^\top\!\! \left(\!\dot{Z} - \begin{bsmallmatrix}
                                \gamma I_n \\ 0
                            \end{bsmallmatrix}E\right)^\top\!\!\! \prec 0\\
                            ZQ = Q^\top Z^\top \succ 0.
                        \end{cases}
                    \end{equation}
                    \State \emph{Control gain:}
                    \begin{equation}\label{eq:K_state}
                        K = UQ (ZQ)^{-1}.
                    \end{equation}
                    \State \hspace{-0.47cm} \textbf{Control Deployment}
                    \State \emph{Control law:}
                    \begin{equation}\label{eq:controller_state}
                        \dot{\hat{\zeta}}_{\text{c}} = -\lambda \hat{\zeta}_{\text{c}} + \gamma \begin{bmatrix}
                            x \\ u
                        \end{bmatrix}, \qquad  u = K\hat{\zeta}_{\text{c}}.
                    \end{equation}
                    Initialization: $\hat{\zeta}_{\text{c}}(0) \in \R^{n+m}$ arbitrary.
        		\end{algorithmic}
        	\end{algorithm}

            The proposed procedure, described in Algorithm \ref{alg:input_state_data}, is based on the following key ideas:
            \begin{itemize}
                \item Consider an input-state trajectory \eqref{eq:dataset_state} of system \eqref{eq:plant_state}.
                Choose gains $\lambda$ and $\gamma$ such that, in addition to $\gamma \neq 0$, $\lambda > 0$.
                By Lemma \ref{lemma:i_o_equivalence_state}, data \eqref{eq:dataset_state} can be seen as an input-output trajectory of system \eqref{eq:filter_state_compact} with $\xi(t) = x(t)$.
                \item Since \eqref{eq:filter_state_compact} is equivalent to \eqref{eq:filter_state}, its behavior is simulated with \eqref{eq:filter_state_obs}, which is a low-pass filter of the data due to $\lambda > 0$ and can be interpreted as a state observer of \eqref{eq:filter_state_compact}. 
                \item The simulated non-minimal state $\hat{\zeta}$ and its derivative $\dot{\hat{\zeta}}$ can be used for a data-driven control strategy that exploits realization \eqref{eq:filter_state_compact} to stabilize the original plant \eqref{eq:plant_state}.
                However, since $A$ and $B$ are unknown, $\hat{\zeta}(0)$ cannot be chosen to perfectly match the trajectories of \eqref{eq:filter_state} and \eqref{eq:filter_state_obs}.
                Therefore, the algorithm needs to account for the fact that \eqref{eq:filter_state_obs} converges only asymptotically to a non-minimal realization of \eqref{eq:plant_state}.
                For this, we choose the special initialization $\hat{\zeta}(0) = 0$.
            \end{itemize}
            We now make the previous arguments precise.
            Consider the interconnection of plant \eqref{eq:plant_state} and filter \eqref{eq:filter_state_obs}, having states $(x, \hat{\zeta})$.
            To characterize the filter transient, define the mismatch error
            \begin{equation}\label{eq:epsilon_state}
                \epsilon \coloneqq x - \frac{1}{\gamma}\begin{bmatrix}
                    A + \lambda I_n & B
                \end{bmatrix}\hat{\zeta},
            \end{equation}
            which originates from the fact that we cannot ensure $x(0) = \gamma^{-1}\begin{bmatrix}A + \lambda I_n & B\end{bmatrix}\hat{\zeta}(0)$.
            The evolution of $\epsilon$ can be computed from \eqref{eq:plant_state}, \eqref{eq:filter_state_obs} as follows:
            \begin{equation}
                \begin{split}
                    \dot{\epsilon} &= Ax + Bu - \frac{1}{\gamma}\begin{bmatrix}
                    A + I_n & B
                \end{bmatrix}\left(-\lambda \hat{\zeta} + \gamma\begin{bmatrix}
                    x \\ u
                \end{bmatrix} \right)\\
                    &= -\lambda x + \frac{\lambda}{\gamma} \begin{bmatrix}
                    A + I_n & B
                \end{bmatrix} \hat{\zeta} = -\lambda \epsilon.
                \end{split} 
            \end{equation}
            Using the change of coordinates \eqref{eq:epsilon_state}, the interconnection of \eqref{eq:plant_state} and \eqref{eq:filter_state_obs} can be represented with states $(\epsilon, \hat{\zeta})$ as:
            \begin{equation}\label{eq:open_loop_state}
                \begin{split}
                    \dot{\epsilon} &= -\lambda \epsilon\\
                    \dot{\hat{\zeta}} &= \underbrace{\begin{bmatrix}
                        A & B\\
                        0 & -\lambda I_m
                    \end{bmatrix}}_{\eqqcolon F}\hat{\zeta} + \underbrace{\begin{bmatrix}
                        0 \\ \gamma I_m
                    \end{bmatrix}}_{\eqqcolon G}u + \underbrace{\begin{bmatrix}
                        \gamma I_{n} \\ 0
                    \end{bmatrix}}_{\eqqcolon D}\epsilon,
                \end{split}
            \end{equation}
            where the $\hat{\zeta}$-subsystem is a system with the same structure of \eqref{eq:filter_state_compact} and subject to the perturbation $D\epsilon$, which converges to $0$ exponentially.
            
            From \eqref{eq:epsilon_state} and choosing $\hat{\zeta}(0) = 0$, it holds that $\epsilon(0) = x(0)$.
            Thus, $\epsilon(t) = e^{-\lambda t}x(0)$ can be computed for every $t \in [0, T]$.
            The proposed procedure involves collecting $N$ samples of $u$, $\hat{\zeta}$, $\dot{\hat{\zeta}}$, and $\epsilon$ as shown in \eqref{eq:batches_state}, then solving LMI \eqref{eq:LMI_state} and computing a control gain $K$ from \eqref{eq:K_state}.
            The resulting controller \eqref{eq:controller_state} is a dynamic feedback law that incorporates the filter dynamics.
            Note that the state $\hat{\zeta}_{\text{c}}$ of \eqref{eq:controller_state} can be initialized arbitrarily.
            We are ready to present the theoretical guarantees for Algorithm \ref{alg:input_state_data}.
            
            \begin{theorem}\label{thm:state}
                Consider Algorithm \ref{alg:input_state_data} and let Assumption \ref{hyp:ctrb} hold.
                Then:
                \begin{enumerate}
                    \item LMI \eqref{eq:LMI_state} is feasible if the batch \eqref{eq:batches_state} is exciting, i.e.:
                    \begin{equation}\label{eq:ZU_rank}
                        \rank\begin{bmatrix}
                            Z \\ U
                        \end{bmatrix} = n + 2m.
                    \end{equation}
                    \item For any solution $Q$ of \eqref{eq:LMI_state}, the gain $K$ computed from \eqref{eq:K_state} is such that  $F + GK$ is Hurwitz.
                    As a consequence, the origin $(x, \hat{\zeta}_{\textup{c}}) = 0$ of the closed-loop interconnection of plant \eqref{eq:plant_state} and controller \eqref{eq:controller_state} is globally exponentially stable.
                \end{enumerate}
            \end{theorem}
            
            \begin{remark}\label{rem:excitation_state}
                Due to space limitations, we do not formally study how to guarantee condition \eqref{eq:ZU_rank}.
                However, we provide some useful insights:
                \begin{itemize}
                    \item System \eqref{eq:filter_state_compact} is controllable by Lemma \ref{lemma:ctrb_state}, so $(\zeta(t), u(t))$ is persistently exciting by \cite[Thms. 1 and 2]{nordstrom1987persistency} if the input $u(t)$ is sufficiently rich.
                    Since $\hat{\zeta}(t) - \zeta(t) \to 0$ exponentially, for a dataset length $T > 0$ sufficiently high, there exists $\mu > 0$ such that:
                    \begin{equation}\label{eq:E_continuous_time}
                        \int_{0}^{T}\begin{bmatrix}
                            \hat{\zeta}(\tau)\\
                            u(\tau)
                        \end{bmatrix}
                        \begin{bmatrix}
                            \hat{\zeta}(\tau)\\
                            u(\tau)
                        \end{bmatrix}^\top \textup{d}\tau \succeq \mu I_{n + 2m}.
                    \end{equation}
                    \item Under sufficient smoothness of the involved signals and sufficiently small sampling time $\Ts > 0$ (see, e.g., \cite[Lemma IV.3]{eising2024sampling}),  \eqref{eq:E_continuous_time}
                    implies \eqref{eq:ZU_rank}.
                \end{itemize}
            \end{remark}

    \section{Data-Driven Control from Input-Output Data}\label{sec:output_feedback}
        To highlight the parallelism with the state-feedback scenario, we slightly abuse the notation of Section \ref{sec:state_feedback} by adopting here similar symbols.
        Consider a single-input single-output system of the form
        \begin{equation}\label{eq:plant_output}
        	\begin{split}
        	    \dot{x} &= A x + b u\\
                y &= c^\top x,
        	\end{split}
        \end{equation}
        where $x \in \R^n$ is the unmeasured state, $u \in \R$ is the control input, $y \in \R$ is the measured output, and $A$, $b$, and $c$ are matrices of appropriate dimensions whose values are unknown but satisfy the following assumption.
        \begin{assumption}\label{hyp:ctrb_obs}
            The pair $(A, b)$ is controllable and the pair $(c^\top, A)$ is observable.
        \end{assumption}
        Compared to Section \ref{sec:state_feedback}, the additional challenge of this scenario is that only the output $y$ is available instead of the state $x$.
        However, once again, it is possible to introduce a filter of the data that reconstructs a non-minimal realization of the plant, thus enabling the application of data-driven control without measuring derivatives.
            
        \subsection{Non-Minimal System Realization}
            Define the following dynamical system, having input $u$ and output $\eta \in \R$:
            \begin{equation}\label{eq:filter_output_compact}
                \begin{split}
                    \dot{\zeta} &= \begin{bmatrix}
                        \Lambda + \ell \theta_1^\top & \ell \theta_2^\top\\
                        0 & \Lambda
                    \end{bmatrix} \zeta + \begin{bmatrix}
                        0 \\ \ell
                    \end{bmatrix} u\\
                    \eta &= \begin{bmatrix}
                        \theta_1^\top & \theta_2^\top
                    \end{bmatrix}\zeta,
                \end{split}
            \end{equation}
            where $\zeta \in \mathbb{R}^{2n}$ is the state , $\Lambda \in \R^{n \times n}$ and $\ell \in \R^n$ are constant tuning gains, and $\theta_1$, $\theta_2 \in \R^n$ are constant vectors whose values, unknown for design, will be derived from $A$, $b$, $c$, $\Lambda$, and $\ell$ to match the input-output behavior of systems \eqref{eq:plant_output} and \eqref{eq:filter_output_compact}.
            System \eqref{eq:filter_output_compact} is obtained from the following representation, also used in the literature of adaptive observers \cite[Ch. 4]{narendra2012stable}:
            \begin{equation}\label{eq:filter_output}
                \begin{split}
                    \dot{\zeta} &=  \begin{bmatrix}
                        \Lambda & 0\\
                        0 & \Lambda
                    \end{bmatrix}\zeta + \begin{bmatrix}
                        \ell & 0\\
                        0 & \ell
                    \end{bmatrix}\begin{bmatrix}
                        \eta \\ u
                    \end{bmatrix}\\
                    \eta &= \begin{bmatrix}
                        \theta_1^\top & \theta_2^\top
                    \end{bmatrix}\zeta   
                \end{split}
            \end{equation}
            which can be interpreted as a filter of $\eta$ and $u$.

            In order to extend the properties found in Section \ref{sec:filter_state} to the input-output scenario, vectors $\theta_1$ and $\theta_2$ can be chosen according to the following fundamental result.
            
            \begin{lemma}\label{lemma:Pi}
                Under Assumption \ref{hyp:ctrb_obs}, for all controllable pairs $(\Lambda, \ell)$ such that $\Lambda$ has distinct eigenvalues, there exist a full rank matrix $\Pi \in \R^{n \times 2n}$ and vectors $\theta_1$, $\theta_2 \in \R^n$ such that:
                \begin{equation}\label{eq:Pi}
                    \begin{split}
                        \Pi \begin{bmatrix}
                            \Lambda + \ell \theta_1^\top & \ell \theta_2^\top\\
                            0 & \Lambda
                        \end{bmatrix} &= A \Pi, \quad \Pi \begin{bmatrix}
                            0 \\ \ell
                        \end{bmatrix} = b\\
                        \begin{bmatrix}
                            \theta_1^\top & \theta_2^\top
                        \end{bmatrix} &= c^\top \Pi.
                    \end{split}
                \end{equation}
                In particular, $A - \Pi_1\ell c^\top$ and $\Lambda$ are similar, where $\Pi_1$ is the matrix given by the first $n$ columns of $\Pi$.
            \end{lemma}

            The next results are equivalent to Lemmas \ref{lemma:i_o_equivalence_state} and \ref{lemma:ctrb_state}.
            \begin{lemma}\label{lemma:i_o_equivalence_output}
                Given the assumptions and matrices $\Pi$, $\theta_1$, $\theta_2$ of Lemma \ref{lemma:Pi}, the observable and controllable subsystem of \eqref{eq:filter_output_compact} obeys the same dynamics of \eqref{eq:plant_output}, with state $\Pi \zeta \in \R^n$ and output $\eta \in \R$.
            \end{lemma}
            \begin{lemma}\label{lemma:ctrb_output}
                Given the assumptions and matrices $\Pi$, $\theta_1$, $\theta_2$ of Lemma \ref{lemma:Pi}, pair
                \begin{equation}\label{eq:AB_filter_output}
                    \left(\begin{bmatrix}
                        \Lambda + \ell \theta_1^\top & \ell \theta_2^\top\\
                        0 & \Lambda
                    \end{bmatrix}, \begin{bmatrix}
                            0 \\ \ell
                        \end{bmatrix} \right)
                \end{equation}
                is controllable.
            \end{lemma}
        
        \subsection{Controller Design}
            \begin{algorithm}[t!]
        		\caption{Controller Design from Input-Output Data}\label{alg:input_output_data}
        		\begin{algorithmic}
        			\State \hspace{-0.47cm} \textbf{Initialization} 
        			\State \emph{Dataset:}
                    \begin{equation}\label{eq:dataset_output}
                        (u(t), y(t)), \qquad \forall t \in [0, T].
                    \end{equation}
                    \State \emph{Tuning:} $\Lambda = \diag(-\lambda_1, \ldots, -\lambda_n )$, with $0 < \lambda_1 < \ldots < \lambda_n$, $\ell = \col(\gamma_1, \ldots, \gamma_n)$, with $\gamma_1, \ldots, \gamma_n \neq 0$, $\Ts > 0$.
        			\State \hspace{-0.47cm} \textbf{Data Batches Construction} 
                    \State \emph{Filter of the data:}
                    \begin{equation}\label{eq:filter_output_obs}
                        \dot{\hat{\zeta}}(t) = \! \begin{bmatrix}
                            \Lambda & 0\\
                            0 & \Lambda
                        \end{bmatrix}\!\hat{\zeta}(t) + \begin{bmatrix}
                            \ell & 0\\
                            0 & \ell
                        \end{bmatrix}\!\!\begin{bmatrix}
                            y(t) \\ u(t)
                        \end{bmatrix}\!,\quad \forall t \in [0, T].
                    \end{equation}
                    Initialization: $\hat{\zeta}(0) = 0$.
                    \State \emph{Auxiliary dynamics:}
                    \begin{equation}\label{eq:chi}
                        \dot{\chi}(t) = \Lambda \chi(t), \qquad \forall t \in [0, T].
                    \end{equation}
                    Initialization: $\chi(0) = [1\; \cdots \; 1]^\top$.
                    \State \emph{Sampled data batches:}
        			\begin{equation}\label{eq:batches_output}
                        \begin{split}
                            U \! &\coloneqq \! \begin{bmatrix}
                            u(0)\! & u(\Ts)\! & \!\cdots\! & u((N \! - \! 1)\Ts)
                            \end{bmatrix} \! \in \! \R^{1 \times N}\\
                            Z_{\text{a}} \! &\coloneqq \! \begin{bmatrix}
                            \begin{bmatrix}
                                \chi(0)\\
                                \hat{\zeta}(0)
                            \end{bmatrix} & \cdots & \begin{bmatrix}
                                \chi((N - 1)\Ts)\\
                                \hat{\zeta}((N - 1)\Ts)
                            \end{bmatrix}
                            \end{bmatrix} \! \in \! \R^{3n\times N} \\
                            \dot{Z}_{\text{a}} \! &\coloneqq \! \begin{bmatrix}
                            \begin{bmatrix}
                                \dot{\chi}(0)\\
                                \dot{\hat{\zeta}}(0)
                            \end{bmatrix} & \cdots & \begin{bmatrix}
                                \dot{\chi}((N - 1)\Ts)\\
                                \dot{\hat{\zeta}}((N - 1)\Ts)
                            \end{bmatrix}
                            \end{bmatrix} \! \in \! \R^{3n \times N}.
                        \end{split}
                    \end{equation}
                    \State \hspace{-0.47cm} \textbf{Stabilizing Gain Computation}
                    \State \emph{LMI:} find $Q \in \R^{N\times 3n}$ such that:
                    \begin{equation}\label{eq:LMI_output}
                         \begin{cases}
                            \dot{Z}_{\text{a}} Q + Q^\top \dot{Z}_{\text{a}}^\top \prec 0\\
                            Z_{\text{a}} Q = Q^\top Z_{\text{a}}^\top \succ 0.
                        \end{cases}
                    \end{equation}
                    \State \emph{Control gain:}
                    \begin{equation}\label{eq:K_output}
                        K = UQ (Z_{\text{a}}Q)^{-1}\begin{bmatrix}
                            0 \\ I_{2n}
                        \end{bmatrix}.
                    \end{equation}
                    \State \hspace{-0.47cm} \textbf{Control Deployment}
                    \State \emph{Control law:}
                    \begin{equation}\label{eq:controller_output}
                        \dot{\hat{\zeta}}_{\text{c}} = \! \begin{bmatrix}
                            \Lambda & 0\\
                            0 & \Lambda
                        \end{bmatrix}\!\hat{\zeta}_{\text{c}} + \begin{bmatrix}
                            \ell & 0\\
                            0 & \ell
                        \end{bmatrix}\!\!\begin{bmatrix}
                            y \\ u
                        \end{bmatrix}\!, \quad  u = K \hat{\zeta}_{\text{c}}.
                    \end{equation}
                    Initialization: $\hat{\zeta}_{\text{c}}(0) \in \R^{2n}$ arbitrary.
        		\end{algorithmic}
        	\end{algorithm}

            The procedure presented in Algorithm \ref{alg:input_output_data} follows similar ideas to those illustrated in the state-feedback scenario.
            Given an input-output trajectory \eqref{eq:dataset_output} of \eqref{eq:plant_output}, define system \eqref{eq:filter_output_obs}, replicating dynamics \eqref{eq:filter_output} with $\eta(t) = y(t)$.
            In \eqref{eq:filter_output_obs}, we choose $\Lambda$ diagonal with distinct negative diagonal entries and $\ell$ with all non-zero entries.
            As a consequence, $\Lambda$ is Hurwitz and thus \eqref{eq:filter_output_obs} is a low-pass filter of the input-output data \eqref{eq:dataset_output}.
            Also, $\Lambda$ has distinct eigenvalues and it can be verified with the PBH test that pair $(\Lambda, \ell)$ is controllable.
            Thus, Lemmas \ref{lemma:Pi}, \ref{lemma:i_o_equivalence_output}, and \ref{lemma:ctrb_output} hold.
            The proposed design can be derived for more general choices of $\Lambda$ and $\ell$, although at the expense of increased notational burden.
            
            Consider $\Pi$, $\theta_1$, $\theta_2$ from Lemma \ref{lemma:Pi}, then define:
            \begin{equation}\label{eq:epsilon_output}
                \epsilon \coloneqq x - \Pi \hat{\zeta},
            \end{equation}
            and note that $\epsilon(0) = x(0)$ since we choose $\hat{\zeta}(0) = 0$.
            The dynamics of $\epsilon$ are computed from \eqref{eq:plant_output}, \eqref{eq:Pi}, \eqref{eq:filter_output_obs}, and \eqref{eq:epsilon_output} as follows:
            \begin{equation}
                \begin{split}
                    \dot{\epsilon} &= Ax + bu - \Pi \begin{bmatrix}
                        \Lambda & 0\\
                        0 & \Lambda
                    \end{bmatrix} \hat{\zeta} - \Pi\begin{bmatrix}
                        \ell & 0\\
                        0 & \ell
                    \end{bmatrix}\begin{bmatrix}
                        c^\top x \\ u
                    \end{bmatrix}\\
                    &= (A - \Pi_1 \ell c^\top)\epsilon + \left(\!A\Pi - \Pi \begin{bmatrix}
                        \Lambda + \ell \theta_1^\top & \ell \theta_2^\top\\
                        0 & \Lambda
                    \end{bmatrix}\!\right)\!\hat{\zeta}\\
                    & = H \Lambda H^{-1} \epsilon,
                \end{split}
            \end{equation}
            where $H$ is a non-singular matrix that exists due to $\Lambda$ and $A - \Pi_1 \ell c^\top$ being similar by Lemma \ref{lemma:Pi}.
            We can write the interconnection of plant \eqref{eq:plant_output} and filters \eqref{eq:filter_output_obs} using the change of coordinates \eqref{eq:epsilon_output}, leading to
            \begin{equation}\label{eq:open_loop_output}
                \begin{split}
                    \dot{\epsilon} &= H \Lambda H^{-1} \epsilon\\
                    \dot{\hat{\zeta}} &= \underbrace{\begin{bmatrix}
                        \Lambda + \ell\theta_1^\top & \ell\theta_2^\top\\
                        0 & \Lambda
                    \end{bmatrix}}_{\eqqcolon F}\hat{\zeta} + \underbrace{\begin{bmatrix}
                        0 \\ \ell
                    \end{bmatrix}}_{\eqqcolon g}u + \underbrace{\begin{bmatrix}
                        \ell c^\top \\ 0
                    \end{bmatrix}}_{\eqqcolon D}\epsilon,
                \end{split}
            \end{equation}
            which shares the same structure of \eqref{eq:open_loop_state}.
            
            Contrary to Section \ref{sec:state_feedback}, $D\epsilon$ is not available in the output-feedback scenario.
            Since $\epsilon \to 0$ exponentially, a simple approach would be to sample $u$, $\hat{\zeta}$, and $\dot{\hat{\zeta}}$ after a sufficiently long time to make the perturbation $D\epsilon$ small enough.
            This method, however, would cause an inefficient use of data and would not be rigorous as $\epsilon$ can be arbitrarily large due to $A$, $b$, $c$, and $x(0)$ being unknown.
            In the following, instead, we propose an approach that compensates $D\epsilon$ exactly without any need for a waiting time.
            
            From \eqref{eq:open_loop_output} and $\epsilon(0) = x(0)$, $\epsilon(t)$ can be computed as
            \begin{equation}
                \epsilon(t) = e^{H\Lambda H^{-1}t}\epsilon(0) = He^{\Lambda t}H^{-1}x(0) = L \chi(t),
            \end{equation}
            where $L \in \R^{n \times n}$ is an unknown matrix depending on $H$ and $x(0)$\footnote{$L = ((H^{-1}x(0))^\top\otimes H)\diag(\mathrm{e}_1, \ldots, \mathrm{e}_n)$, where $\otimes$ denotes the Kronecker product and $(\mathrm{e}_1, \ldots, \mathrm{e}_n)$ is the orthonormal basis of $\R^n$.}, while $\chi(t) \coloneqq [e^{-\lambda_1 t}\; \cdots \; e^{-\lambda_n t}]^\top \in \R^n$.
            Note that $\chi(t)$ obeys dynamics \eqref{eq:chi}, with $\chi(0) = [1\; \cdots \; 1]^\top$.
            Thus, the sequence $(u(t), \hat{\zeta}(t))$ obtained from \eqref{eq:dataset_output},  \eqref{eq:filter_output_obs} satisfies for all $t \in [0, T]$ the following differential equation:
            \begin{equation}\label{eq:virtual_system}
                \begin{bmatrix}
                    \dot{\chi}\\
                    \dot{\hat{\zeta}}
                \end{bmatrix} = \begin{bmatrix}
                    \Lambda & 0\\
                    DL      & F
                \end{bmatrix}\begin{bmatrix}
                    \chi \\ \hat{\zeta}
                \end{bmatrix} + \begin{bmatrix}
                    0 \\ g
                \end{bmatrix}u,
            \end{equation}
            with initial conditions $\chi(0) = [1\; \cdots \; 1]^\top$, $\hat{\zeta}(0) = 0$.

            Note that the state and the state derivative of \eqref{eq:virtual_system} are available.
            As a consequence, we can sample $u$,  $(\chi, \hat{\zeta})$, and $(\dot{\chi}, \dot{\hat{\zeta}})$, to obtain the batches \eqref{eq:batches_output}.
            Then, LMI \eqref{eq:LMI_output} can be used to compute a feedback law for plant \eqref{eq:plant_output}.
            Suppose that there exists a matrix $Q$ that solves LMI \eqref{eq:LMI_output}, then we can use the same result of Theorem \ref{thm:state}, point 2, to ensure that $[K_\chi \; K] = UQ(Z_{\text{a}}Q)^{-1}$ is such that
            \begin{equation}
                \begin{bmatrix}
                    \Lambda & 0\\
                    DL + gK_\chi & F + gK 
                \end{bmatrix}
            \end{equation}
            is Hurwitz.
            Thus, $K$ computed from \eqref{eq:K_output} is such that $F + gK$ is Hurwitz.
            The resulting observer-based controller is given in \eqref{eq:controller_output}.
            
            Given the steps presented in the discussion above, we are ready to state the theoretical result for Algorithm \ref{alg:input_output_data}.
            
            \begin{theorem}\label{thm:output}
            Consider Algorithm \ref{alg:input_output_data} and let Assumption \ref{hyp:ctrb_obs} hold.
            Then:
                \begin{enumerate}
                    \item LMI \eqref{eq:LMI_output} is feasible if the batch \eqref{eq:batches_output} is exciting, i.e.:
                    \begin{equation}\label{eq:ZaU_rank}
                        \rank\begin{bmatrix}
                            Z_{\text{a}} \\ U
                        \end{bmatrix} = 3n + 1.
                    \end{equation}
                    \item For any solution $Q$ of \eqref{eq:LMI_output}, the gain $K$ computed from \eqref{eq:K_output} is such that  $F + gK$ is Hurwitz.
                    As a consequence, the origin $(x, \hat{\zeta}_{\textup{c}}) = 0$ of the closed-loop interconnection of plant \eqref{eq:plant_output} and controller \eqref{eq:controller_output} is globally exponentially stable.
                \end{enumerate}
            \end{theorem}

            \begin{remark}
                $V \coloneqq [\chi(0)\; \cdots \; \chi((N-1)\Ts)]$ is a Vandermonde matrix with roots $e^{-\lambda_1 \Ts}, \ldots, e^{-\lambda_n \Ts}$, so it has full row rank when $N \geq n$.
                Define $Z \coloneqq [\hat{\zeta}(0)\; \cdots \; \hat{\zeta}((N-1)\Ts)]$.
                Then, \eqref{eq:ZaU_rank} holds if $[Z^\top\, U^\top]^\top$ has full row rank and is linearly independent from $V$.
                For the full rank of $[Z^\top\, U^\top]^\top$, we refer to Remark \ref{rem:excitation_state}.
                To give an intuition on the second requirement, let $\hat{\zeta}(t)$ and $u(t)$ be the sum of $p$ distinct sinusoids.
                From $\sin(\omega t) = (e^{\mathrm{i}\omega t} - e^{-\mathrm{i}\omega t})/(2\mathrm{i})$, $\cos(\omega t) = (e^{\mathrm{i}\omega t} + e^{-\mathrm{i}\omega t})/2$:
                \begin{equation}
                    \begin{bmatrix}
                        Z \\ U
                    \end{bmatrix}\! = \! \Psi W = \! \Psi \! \begin{bmatrix}
                        1 & e^{\mathrm{i} \omega_1 \Ts} & \cdots &  e^{\mathrm{i} \omega_1 (N-1) \Ts} \\
                        1 & e^{-\mathrm{i} \omega_1 \Ts} & \cdots &  e^{-\mathrm{i} \omega_1 (N-1) \Ts} \\
                        \vdots & \vdots & & \vdots \\
                        1 & e^{\mathrm{i} \omega_{p} \Ts} & \cdots &  e^{\mathrm{i} \omega_{p} (N-1) \Ts} \\
                        1 & e^{-\mathrm{i} \omega_{p} \Ts} & \cdots &  e^{-\mathrm{i} \omega_{p} (N-1) \Ts} \\
                    \end{bmatrix},
                \end{equation}
                where $\Psi \in \mathbb{C}^{(2n + 1)\times 2p}$ and $W$ is a Vandermonde matrix.
                As a consequence, $N \geq n + 2p$ ensures that $W$ and $V$ are linearly independent, implying that each non-zero row of $\Psi W$ is linearly independent from $V$.
            \end{remark}

    \section{Numerical Examples}\label{sec:numerical_example}
        Algorithms \ref{alg:input_state_data} and \ref{alg:input_output_data} have been implemented in MATLAB using YALMIP \cite{Lofberg2004} and MOSEK \cite{mosek} to solve the LMIs.
        The developed code is available at the linked repository\footnote{\url{https://github.com/IMPACT4Mech/continuous-time_data-driven_control}}.
        
        \subsection{Design with Input-State Data}
            We consider the continuous-time linearized model of an unstable batch reactor given in \cite{walsh2001scheduling}, also used in \cite{de2019formulas} after time discretization.
            The system matrices are
            \begin{equation}
                A = 
                \begin{bmatrix*}[r]
                    1.38& -0.2077& 6.715& -5.676\\
                    -0.5814& -4.29& 0& 0.675\\
                    1.067& 4.273& -6.654& 5.893\\
                    0.048& 4.273& 1.343& -2.104
                \end{bmatrix*}, \quad
                B =
                \begin{bmatrix*}[r]
                    0& 0\\
                    5.679& 0\\
                    1.136& -3.146\\
                    1.136& 0
                \end{bmatrix*},
            \end{equation}
            where $(A, B)$ is controllable and the eigenvalues of $A$ are $\{-8.67, -5.06, 0.0635, 1.99\}$.
            We consider an exploration interval of length $T = 1.5$ s, where we apply a sum of $4$ sinusoids to both entries of $u$ and select $8$ distinct frequencies to ensure informative data.
            We choose filter gains $\lambda = \gamma = 1$ and collect the data with sampling time $\Ts = 0.1$ s.

            Algorithm \ref{alg:input_state_data} has been extensively tested for random initial conditions $x(0)$, with each entry extracted from the uniform distribution $\mathcal{U}(-1, 1)$, returning each time a stabilizing controller.
            For the case $x(0) = [0.311\;\; {-0.6576}\;\; 0.4121\;\; {-0.9363}]^\top$, we obtain the gain
            \begin{equation}
                K =
                \begin{bmatrix*}[r]
                    -1.507&  -18.69&    0.155&   -0.681&    2.925&    0.79\\
                    17.45&    0.224&   44.06&  -36.37&   1.09&   -3.518
                \end{bmatrix*}\!\!,
            \end{equation}
            which places the eigenvalues of the closed-loop matrix $F + GK$ in $\{-5.107 \pm 10.729\mathrm{i}, -1.238, -1.024 \pm 9.654\mathrm{i}, -0.759\}$.
        
        \subsection{Design with Input-Output Data}
            We consider a non-minimum phase SISO system having input-output behavior specified by the transfer function 
            \begin{equation}
                c^\top(sI - A)^{-1}b = \frac{s-1}{s(s^2 + 4)},
            \end{equation}
            for which we choose a minimal realization in controllability canonical form.
            We perform an exploration of $T = 2$ s with an input $u$ given by the sum of $4$ sinusoids at distinct frequencies.
            We choose filter parameters $\Lambda = \diag(-1, -2, -3)$, $\ell = \col(1, 2, 3)$ and sampling time $\Ts = 0.1$ s.

            Similar to the previous case, Algorithm \ref{alg:input_output_data} has been extensively tested with random initial conditions so that each entry of $x(0)$ is extracted from the uniform distribution $\mathcal{U}(-5, 5)$.
            In each test, the procedure returned a stabilizing controller.
            For the case $x(0) = [{-3.9223}\;\; 4.0631\;\; 3.7965]^\top$, we obtain the gain
            \begin{equation}
                K =
                \begin{bmatrix}
                    -0.508&    3.208&   -2.392&    0.001&   -0.577&    1.055
                \end{bmatrix},
            \end{equation}
            which places the eigenvalues of $F + gK$ in $\{-2.028, -0.723 \pm 0.647\mathrm{i}, -0.22, -0.147 \pm 2.09\mathrm{i}\}$.   

    \section{Conclusion}\label{sec:conclusion}
        We addressed the problem of data-driven stabilization of unknown continuous-time linear time-invariant systems by proposing a framework that combines signal filtering with LMIs.
        Specifically, we employed low-pass filters that reconstruct a non-minimal realization of the plant.
        Then, LMIs inspired by those of \cite{de2019formulas} and based on the non-minimal realization have been used to compute the gains of a dynamic, filter-based, controller.
        This approach circumvents the need for signal derivatives without resorting to noise-sensitive numerical techniques like finite differences.
        Future research will focus on extending the method to multi-input multi-output and nonlinear systems, as well as exploring the conditions to ensure exciting data.

    \bibliographystyle{IEEEtran}
	
	\bibliography{data-driven_bib}

\end{document}